\documentclass[a4paper,12pt]{amsart}

\usepackage{amsmath}
\usepackage{amssymb}
\usepackage{mathrsfs}
\usepackage{enumerate}
\usepackage{ifthen}
\usepackage{graphicx}
\usepackage[T1]{fontenc} 

\nonstopmode \numberwithin{equation}{section}
\newtheorem{theorem}{Theorem}[section]
\newtheorem{corollary}{Corollary}[section]
\newtheorem{lemma}{Lemma}[section]
\newtheorem{prop}[equation]{Proposition}
\theoremstyle{definition}

\newtheorem{problem}{Problem}[section]

\begin{document}

\title{On distance magic harary graphs}

\author{A V Prajeesh}
\address{A V Prajeesh, Department of Mathematics, National Institute of Technology Calicut, Kozhikode~\textnormal{673601}, India. }
\email{prajeesh\_p150078ma@nitc.ac.in}

\author{K Paramasivam}
\address{Krishnan Paramasivam, Department of Mathematics, National Institute of Technology Calicut,  Kozhikode~\textnormal{673601}, India. }
\email{sivam@nitc.ac.in}



\begin{abstract}
	This paper establishes two techniques to construct larger distance magic and $(a,d)$-distance antimagic graphs using Harary graphs and  provides a solution to the existence of distance magicness of $G\circ C_4$ and $G\times C_4$, for every non-regular distance magic graph $G$  with maximum degree $|V(G)|-1$.
\end{abstract}

\thanks{}

\maketitle
\pagestyle{myheadings}
\markboth{A V Prajeesh and Krishnan Paramasivam }{On Distance Magic Harary Graphs}

\section{Introduction}
\noindent In this paper, we consider only simple and finite graphs. We use $V(G)$ for the vertex set and $E(G)$ for the edge set of a graph $G$. The neighborhood $N_G(v)$, or shortly $N(v)$ of a vertex $v$ of $G$ is the set of all vertices adjacent to $v$. For standard graph theoretic notations and definitions, we refer Bondy and Murty \cite{bondy1976graph}, and Hammack $et$ $al.$\cite{hammack2011handbook}.
\par A distance magic labeling of $G$ is a bijection $f : V(G) \rightarrow \{1,2,...,|V(G)|\}$, such that for any $u$ of $G$, the weight of $u$, $w_G(u) = \sum\limits_{v\in N_{G}(u)} f(v) $ is a constant $c$. A graph $G$ that admits such a labeling is called a distance magic graph.
\par The concept of distance magic labeling was studied by Vilfred \cite{vilfred} as sigma labeling. Later, Miller $et$ $al.$ \cite{miller2003distance} called it a $1$-vertex magic vertex labeling and Sugeng $et$ $al.$ \cite{sugeng2009distance} referred the same as distance magic labeling.
\par An equalized incomplete tournament of $n$ teams with $r$ rounds,
$EIT(n, r)$ is a tournament which satisfies the following conditions:
\begin{itemize}
	\item [(i)] every team plays against exactly $r$ opponents.
	\item [(ii)] the total strength of the opponents, against which each team plays is a constant.
\end{itemize}
\par Therefore, finding a solution for an equalized incomplete tournament $EIT(n, r)$
is equivalent to establish a distance magic labeling of an $r$-regular graph of
order $n$. For more details, one can refer \cite{froncek2006fair,froncek2007fair}.
\par The important results and problems, which are relevant and helpful in proving our results, are listed below.
\begin{theorem}\label{oddregular}
	\textnormal{\cite{vilfred,miller2003distance,jinnah,rao}} No $r$-regular graph with $r$-odd can be a distance magic graph.
\end{theorem}
\begin{lemma}\label{degreenecessary}
	\textnormal{\cite{miller2003distance}} If $G$ contains two vertices $u$ and $v$ such that $|N_{G}(u) \cap N_{G}(v)| = d(v)-1 = d(u)-1$, then $G$ is not distance magic.
\end{lemma}
\begin{theorem}\label{completedmg}
	\textnormal{\cite{miller2003distance}} $K_{n}$ is distance magic if and only if $n=1$.
\end{theorem}
\par A distance magic graph $G$ on $n$ vertices, is called balanced if there exists a bijection $f:V(G) \rightarrow \{1,2,...,n\}$ such that for any $w$ of $G$, the following holds: if $u \in N_{G}(w)$ with $f(u)=i$, then there exists $v \in N_{G}(w)$ and $v \neq u$ such that $f(v)= n+1-i$. We call $u$ and $v$ are the twin vertices to each other. Also, the label-sum of twin vertices, $f(u)+f(v)$ is equal to $n+1$.
\par From \cite{anholcer2015distance}, we observe that $G$ is a balanced distance magic graph if and only if $G$ is regular and the vertex set of $G$ can be expressed as $\{v_{i},v_{i}': 1\leq i\leq \frac{n}{2}\}$ such that  for any $i$, $N_G(v_{i}) = N_G(v_{i}')$, where $v_i$ and  $v_i'$ are the twin vertices. The graphs $K_{2n}-M$, $M$ any perfect matching of $K_{2n}$ and $K_{2n,2n}$ are examples of balanced distance magic graphs.\par
The $m^{\textnormal{th}}$ power of a graph $G$ is the graph $G^{m}$ with the same set of vertices as $G$ and any two vertices $u$ and $v$ are connected by an edge if and only if $d_{G}(u,v)\leq m$.\par
In 2016, Cichacz\cite{cichacz2013distance} studied the distance magic labeling of the graph $ C_{n}^{m}$ and proved the following results.
\begin{lemma}\label{lem1}
	\textnormal{\cite{cichacz2013distance}} If $n = 2m + 2$, then $C_{n}^{m}$ is a distance magic graph.
\end{lemma}
\begin{theorem}\label{thm1}
	\textnormal{\cite{cichacz2013distance}} If $m$ is odd, then $C_{n}^{m}$ is a distance magic graph if and only if $2m(m+1) \equiv 0 ~\textnormal{mod}~n,$ $n \geq 2m+2$ and $\frac{n}{gcd(n,m+1)} \equiv 0 ~\textnormal{mod}~2.$
\end{theorem}
\begin{theorem}\label{thm2}
	\textnormal{\cite{cichacz2013distance}} The graph $C_{n}^{2}$ is not a distance magic graph unless $n = 6$.
\end{theorem}
\par Bondy and Murty \cite{bondy1976graph} constructed an $m$-connected graph $H_{m,n}$ on $n$ vertices that has exactly $\lceil \frac{mn}{2} \rceil$ edges. The structure of $H_{m,n}$ depends on the parities of $m$ and $n$; there are three cases.\\
Case 1. If $m$ is even, then $H_{m,n}$ is constructed as follows. It has vertices $0,1,...,n-1$ and two vertices $i$ and $j$ are
joined if $i - \frac{m}{2} \leq j \leq i +  \frac{m}{2}$ (where addition is taken modulo $n$).\\
Case 2. If $m$ is odd and $n$ is even, then $H_{m,n}$ is constructed by first drawing $H_{m-1,n}$ and then adding edges joining vertex $i$ to vertex $i + \frac{n}{2}$ for $1 \leq i \leq \frac{n}{2}$.\\
Case 3. If $m$ and $n$ are odd, then $H_{m,n}$ is constructed by first drawing $H_{m-1,n}$ and then adding edges joining vertex $0$ to vertices $\frac{n-1}{2}$, $\frac{n+1}{2}$ and vertex $i$ to vertex $i + \frac{n+1}{2}$ for $1 \leq i < \frac{n-1}{2}.$ \par In this paper, the vertices $ 0,1,...,n-1 $ of Harary graph, $H_{m,n}$ are renamed as  $v_{0},v_{1},...,v_{n-1}$. Clearly for $m=n-1, H_{m,n}\cong K_{n}$ and for all $m,n$, $H_{2m,n} \cong C_{n}^{m}$. When $m = 2,3$ and $5$, the distance magic labeling of $C_{n}^{m}$ is completely characterized in \cite{cichacz2013distance}. In addition to that Lemma \ref{lem1}, Theorem \ref{thm1} and \ref{thm2} confirms that distance magic Harary graphs always contain a subclass of $C_{n}^{m}$.\par
Shafiq $et$ $al.$\cite{shafiq2009distance} proved the following result and posted a problem on the existence of distance magic labeling of product $G\circ C_{4}$, for any non-regular graph $G$.
\begin{theorem}
	\textnormal{\cite{shafiq2009distance}} Let $r \geq 1$ and $n \geq 3$. If $G$ is an $r$-regular graph and $C_{n}$ the cycle of length $n$, then $G\circ C_{n}$ admits a labeling if and only if $n = 4$.
\end{theorem}
\begin{problem}\label{prob1}
	\textnormal{\cite{shafiq2009distance}} If $G$ is a non-regular graph, determine if there is a distance magic labeling of $G \circ C_{4}$.
\end{problem}
\par Cichacz and G{\"o}rlich \cite{cichacz2018constant}, posted a similar problem in 2018.
\begin{problem}\label{prob2}
	\textnormal{\cite{cichacz2018constant}} If $G$ is non-regular graph, determine if there is a distance magic labeling of $G \times C_{4}$.
\end{problem}
\par Arumugam $et$ $al.$\cite{arumugam2016distance}, discussed the following result as a characterization of entire class of non-regular distance magic graphs $G$ with $\Delta(G) = n-1.$
\begin{theorem}\label{aru1}
	\textnormal{\cite{arumugam2016distance}}
	Let $G$ be any graph of order $n$ with $\Delta(G)= n-1.$ Then $G$ is
	a distance magic graph if and only if $n$ is odd and $G \cong (K_{n-1}-M) + K_{1}$,
	where $M$ is a perfect matching of $K_{n-1}.$
\end{theorem}
\par Later, Arumugam and Kamatchi \cite{arumugam2012d} generalized the concept of distance magic labeling to $(a,d)$-distance antimagic labeling. An $(a,d)$-distance antimagic labeling of a graph $G$ is defined as a bijection $f:V(G)\rightarrow \{1,2,...,n\}$ such that the set of all vertex weights is $\{a,a+d,a+2d,...,a+(n-1)d\}$, where $a$ and $d$ are fixed integers with $d\geq0$. Any graph that admits such a labeling is called an $(a,d)$-distance antimagic graph. Further the condition is relaxed in \cite{kamatchi2013distance} and defined that if $w_{G}(u)\neq w_{G}(v)$, for any two distinct vertices of $G$, then $f$ is called as distance antimagic labeling of $G$. The following problem was posted in \cite{kamatchi2013distance}.
\begin{problem}\label{prob3}
	\textnormal{\cite{kamatchi2013distance}} If $G$ is distance antimagic, is it true that the graphs $G+K_{1}$, $G+K_{2}$ and the Cartesian product $G\Box K_{2}$ are distance antimagic.
\end{problem}
\begin{lemma}\label{siman1}
	\textnormal{\cite{simanjuntak2013distance}} Let $G$ be an $r$-regular graph. If $G$ is $(a,d)$-distance antimagic, then $d\leq r (\frac{n-r}{n-1})$ and $a=\frac{r(n+1)-d(n-1)}{2}.$
\end{lemma}
\par The following terminologies are used in this paper.\par
Hereafter, the graph $(K_{n-1}-M)+K_{1}$, is denoted by $\mathcal{G}=(V(\mathcal{G}),E(\mathcal{G}))$, where $V(\mathcal{G}) = \{v_0,v_1,...,v_{n-1}\}$ and $E(\mathcal{G}) = \{v_{i}v_{j}: 1\leq i,j \leq n-1,i\neq j, j+\frac{n-1}{2}\mod (n-1)\} \cup \{v_{0}v_{i}: 1\leq i\leq n-1\}$. Also, the vertex set of $\mathcal{G}\circ C_{4}$  and $\mathcal{G}\times C_{4}$ are defined as $\cup_{j=0}^{n-1}A_{j}$, where $A_{j} = \{ v_{j}^{i}: i=0,1,2,3\}$ and $v_{j}^{0},v_{j}^{1},v_{j}^{2},v_{j}^{3}$ are the successive vertices of $C_4$ corresponding to the $j^{th}$ vertex of $\mathcal{G}$.\par
To replace a subgraph $H$ of the given graph $G$ by a set $U =\{u_0,u_1,...,u_{k-1}\}$ of $k$ vertices, we mean first remove all the edges of $H$ and then join all the vertices of $U$ to all the vertices of $H$. After the replacement of $H$ by $U$ in $G$, the new graph is denoted by $G^{\dagger}$. Further, the integer-valued functions,
\[ \alpha(i)= \left\{
\begin{array}{ll}
0 &  $ for $ i $  even $ \\
1  &  $ for $ i $ odd $ \\
\end{array}
\right.\]
and \[ \beta_{n}(i)= \left\{
\begin{array}{ll}
0 & $ for $ i<n+1 \\
1  & $ for $ i\geq n+1 \\
\end{array}
\right.\]
are used to define the labeling in a precise way.
\section{Construction of new distance magic graphs using $H_{m,n}$}
\noindent In 2012, Kovar $et$ $al.$\cite{kovar2012note} established a recursive technique to construct a new distance magic graph from an existing $4$-regular distance magic graph. This section gives another useful technique to construct a larger $(2k+2)$-regular distance magic graph $G^{\dagger}$ from an existing graph $G$. The technique mainly invokes few structural properties of balanced distance magic graphs and the regularity is preserved during the entire process. The following Proposition \ref{prop1} and \ref{prop2}, discussed in \cite{cichacz2013distance} and \cite{arumugam2016distance} are used in the construction procedure later in this section.
\begin{prop}\label{prop1}
	$H_{2n,2n+2}$ is distance magic.
\end{prop}
\begin{proof}
	Consider a function $f$,
	\[  f(v_{i}) = \left\{
	\begin{array}{ll}
	2n+2 - i  & $ for $ 0\leq i \leq n \\
	i-n   & $ for $ n+1 \leq i \leq 2n+1. \\
	\end{array}
	\right.\]
	Here, the weight of every vertex of $H_{2n,2n+2}$ is $2n^2+3n$.
\end{proof}
\begin{prop}\label{prop2}
	$H_{2n+1,2n+3}$ is distance magic.
\end{prop}
\begin{proof}
	Consider a function $f$,
	\[  f(v_{i}) = \left\{
	\begin{array}{ll}
	2n+3 - i  & $ for $ 0\leq i \leq n+1 \\
	i-(n+1)   & $ for $ n+2 \leq i \leq 2n+2. \\
	\end{array}
	\right.\]
	Now, the weight of every vertex of $H_{2n+1,2n+3}$ is $2n^2+5n+3$.
\end{proof}
\par Notice that when $m$ is odd and $n$ is even, $H_{m,n}$ is an odd regular graph and hence it is not distance magic. Further, when $m$ and $n$ both are odd, $H_{m,n}$ is a non-regular graph on $n$ vertices. The following theorem characterizes the distance magic labeling of $H_{3,n}$.
\begin{theorem}
	$H_{3,n}$ is distance magic if and only if $n = 5.$
\end{theorem}
\begin{proof}
	$H_{3,5}$ is distance magic by Proposition \ref{prop2}. For all even $n, H_{3,n}$ is not distance magic by Theorem \ref{oddregular} and for all odd $n\neq5$, apply Lemma \ref{degreenecessary} to $H_{3,n}$ by fixing $u = v_{\frac{n-3}{2}}$ and $v = v_{n-2}$.
\end{proof}
\par One can observe that $H_{2n-2,2n}$ is a balanced distance magic graph because $H_{2n-2,2n}$ is isomorphic to $K_{2n}-M$, with $M$ any perfect matching of $K_{2n}$.
\begin{theorem}\label{constthm}
	Let $G$ be a $2n$-regular distance magic graph on $m$ vertices. If $G$ has a subgraph $H_{2n-2,2n}$ such that the label-sum of twin vertices with respect to $H_{2n-2,2n}$ is $m + 1$, then there exists a $2n$-regular distance magic graph $G^{\dagger}$ on $(m+2n-2)$ vertices.
\end{theorem}
\begin{proof}
	Let $H\cong H_{2n-2,2n}$ be a subgraph of a distance magic graph $G$ such that the label-sum of twin vertices with respect to $H$ is $m+1$.
	To construct the graph $G^{\dagger}$, replace $H$ of $G$ by a set $U = \{u_{0},u_{1},...,u_{2n-3}\}$ of $2n-2$ vertices.	
	Consider the function $f^{\dagger}$ defined on $V(G^{\dagger})$ as, \[  f^{\dagger}(x) = \left\{
	\begin{array}{ll}
	n+f(x)-1 & $ for $ x \in V(G) \\
	m\beta_{n-2}(j)+j+1 & $ for $  x = u_{j}$ and $0 \leq j \leq 2n-3. \\
	\end{array}
	\right.\]
	\par For any vertex $x$ of $U$, since $x$ is adjacent to all the vertices of $H$, we get,
	\begin{eqnarray*}
		w_{G^{\dagger}}(x) = \sum\limits_{v\in V(H)}f^{\dagger}(v) =  2n(n-1)+ \sum\limits_{v \in V(H)}f(v) = 2n^2+mn-n.
	\end{eqnarray*}
	\par For any vertex $x$ of $H$, the $(n-1)$-pairs of its twin neighbors $v$ and $v'$ having label-sum $m+1$, are replaced with new $(n-1)$-pairs of vertices $u$ and $u'$ having label-sum $m+2n-1$ and the labels of remaining neighbors of $x$ are increased by $n-1$. Hence,
	\begin{eqnarray*}
		~~~w_{G^{\dagger}}(x) &=& w_G(x)-(n-1)(m+1)+(n-1)(m+2n-1)+2n-2  \\
		&=& w_G(x)+2n^2-2n \\
		&=& 2n^2+mn-n.
	\end{eqnarray*}
	\par Further, for all remaining vertices $x$ of $G^{\dagger}$, the existing weight of $x$ in $G$ is increased by $2n(n-1)$. Hence, $w_{G^{\dagger}}(x) = w_G(x)+2n(n-1)=
	2n^2+mn-n$. Thus $f^{\dagger}$ is a distance magic labeling of $G^{\dagger}$.
\end{proof}
\par Apply the above method to $H_{2n,2n+2}$, one can get the following new class of distance magic graph.
\begin{corollary}
	$H^{\dagger}_{2n,2n+2}$ is distance magic, where $n>1$.
\end{corollary}
\begin{proof}
	The proof follows from the fact that $H_{2n,2n+2} \cong H_{2n-2,2n}+K_{2}^{c}$.
\end{proof}
\par The following theorem provides a subclass of non-regular distance magic graphs by employing the similar technique in Theorem \ref{constthm}. Also, both the existing and newly constructed distance magic graphs are of odd order.
\begin{theorem}
	$H^{\dagger}_{2n+1,2n+3}$ is distance magic, where $n>0$.
\end{theorem}
\begin{proof}
	Let $G \cong H_{2n+1,2n+3}$. Using the fact that, $G \cong H_{2n,2n+2}+K_{1} $, construct the graph $G^{\dagger}$ by replacing the subgraph $H_{2n,2n+2}$ of $G$ by a set $U = \{u_{0},u_{1},...,u_{2n-1}\}$ of $2n$ vertices. If $f$ is the distance magic labeling of $G$ as given in Proposition \ref{prop2}, then the required distance magic labeling $f^{\dagger}$ of $G^{\dagger}$ is defined as,
	\[ f^{\dagger}(x) = \left\{
	\begin{array}{ll}
	2n+f(v_{i})-n\beta_{0}(i) & $ for $ x = v_{i}$ and $0 \leq i \leq 2n+2 \\
	(f(v_{0})-1)\beta_{n-1}(j)+j+1 & $ for $ x = u_{j} $ and $ 0\leq j \leq 2n-1. \\
	\end{array}
	\right.\]
	\par For each $i\in \{1,2,...,2n+2\}$, the vertex $v_{i}$ is adjacent to $v_{0}$ and all the vertices of $\displaystyle\cup_{j=0}^{n-1} \{u_{j},u_{2n-j-1}: f^{\dagger}(u_{j})+f^{\dagger}(u_{2n-j-1})=4n+3 \}$. Also, for each $k\in \{0,1,...,2n-1\}$, both vertices $u_{k}$ and $v_{0}$ are adjacent to all the vertices of $\displaystyle\cup_{l=1}^{n+1} \{v_{l},v_{n+l+1}: f^{\dagger}(v_{l})+f^{\dagger}(v_{n+l+1})=4n+3 \}$. Thus, the weight of any vertex of $G^{\dagger}$ is $4n^{2}+7n+3$.
\end{proof}
\section{Distance magicness of $\mathcal{G}\circ C_{4}$ and $\mathcal{G}\times C_{4}$}
This section discusses Problem \ref{prob1} and \ref{prob2} for all graphs isomorphic to $\mathcal{G}$.
\begin{theorem}\label{prbsol1}
	$\mathcal{G} \circ C_4 $ is not distance magic.
\end{theorem}
\begin{proof}
	On the contrary, assume that $\mathcal{H} \cong \mathcal{G}\circ C_4$ has a distance magic labeling $f$ with magic constant $c$.
	Then, for any $j\in\{1,2,...,n-1\}$, there exist positive integers $a$ and $b_{j}$  such that
	\begin{eqnarray*}
		f( {v}^0_0)+f( {v}^2_0)&=& f( {v}^1_0)+f( {v}^3_0)\hskip .2cm = \hskip .2cm  a, \\
		f( {v}^0_j)+f( {v}^2_j)&=& f( {v}^1_j)+f( {v}^3_j)\hskip .2cm =  \hskip .2cm b_{j}.
	\end{eqnarray*}
	\par Now, for $i\in \{0,1,2,3\}$, the weights of the vertices
	\begin{eqnarray}\label{eqn2}
	w_{\mathcal{H}}( {v}^i_0) &=& a+ 2\sum_{j=1}^{n-1}b_{j}.
	\end{eqnarray}
	\par Also, for $k\in \{1,2,...,n-1\}$, the weights of the vertices
	\begin{eqnarray}\label{eqn3}
	w_{\mathcal{H}}( {v}^i_k) &=& 2a+ b_{k}+2\sum_{j}b_{j},
	\end{eqnarray}
	where $j$ varies from $1$ to $n-1$, $j\neq k$ and $j\neq (k+\frac{n-1}{2}) \mod (n-1)$.\\
	For $j\in\{1,2,...,\frac{n-1}{2}\}$, comparing the equations $w_\mathcal{H}\big( {v}^i_{j}) = w_\mathcal{H}( {v}^i_{j+\frac{n-1}{2}}\big)$, we get
	\begin{equation*}
	b_{j}=b_{j+\frac{n-1}{2}}.
	\end{equation*}
	From \ref{eqn3}, we obtain,
	\begin{eqnarray*}
		w_{\mathcal{H}}\big( {v}^i_{j}) &=&2a+b_{j}+4\sum_{l=1}^{\frac{n-3}{2}} b_{j+l}\\
		w_{\mathcal{H}}\big( {v}^i_{j+\frac{n-3}{2}}) &=&2a+b_{j+\frac{n-3}{2}}+4\sum_{l=1}^{\frac{n-3}{2}} b_{j+\frac{n-3}{2}-l}\\
	\end{eqnarray*}
	and hence, $ b_{j} = b_{j+\frac{n-3}{2}}$, for $j\in \{1,2,...,\frac{n-1}{2}\}$. So assume, $b_{j}=b$ (say) for all $j$. From \ref{eqn2} and \ref{eqn3}, we have
	\begin{equation}\label{eqn4}
	a+2 \sum_{j=1}^{n-1}b_{j} = a+2 \sum_{j=1}^{n-1}b = 2a + (2n-5)b = c\hskip .2cm  \text{or} \hskip .2cm  a = 3b.
	\end{equation}
	\par On the other hand, the sum of the labels of all the vertices, $1+2+...+4n = 2a+2b(n-1)$. By using \ref{eqn4}, $a=\frac{3n(4n+1)}{n+2}$. But the highest possible value for $a$ is $8n-1$. In that case, there is no $n$ satisfying $3n(4n+1) \leq (8n-1)(n+2)$.
\end{proof}
\begin{corollary}
	$H_{2n+1,2n+3} \circ C_{4}$ is not distance magic. \qed
\end{corollary}
\begin{theorem}\label{prbsol2}
	$\mathcal{G}\times C_4 $ is distance magic if and only if $n=5.$
\end{theorem}
\begin{proof}
	Suppose that $\mathcal{H} \cong \mathcal{G}\times C_4$ has a distance magic labeling $f$ with magic constant $c$.
	For $j\in \{1,2,...,n-1\}$, let $s_{1},s_{2},a_{j},b_{j}$  be positive integers such that
	\begin{eqnarray*}
		f( {v}^0_0)+f( {v}^2_0)&= s_{1}  \textnormal{ and} & f( {v}^1_0)+f( {v}^3_0)= s_{2}, \\
		f( {v}^0_j)+f( {v}^2_j)&= a_{j}  \textnormal{ and} & f( {v}^1_j)+f( {v}^3_j) = b_{j}.
	\end{eqnarray*}
	Notice that,
	\begin{eqnarray*}
		w_{\mathcal{H}}( {v}^1_0) &=& w_{\mathcal{H}}( {v}^3_0) = \sum_{j=1}^{n-1}a_{j} ,\\
		w_{\mathcal{H}}( {v}^0_0) &=& w_{\mathcal{H}}( {v}^2_0) = \sum_{j=1}^{n-1}b_{j}.
	\end{eqnarray*}
	For $k\in \{1,2,...,n-1\}$ we have
	\begin{eqnarray*}
		w_{\mathcal{H}}({v}^1_k) &=& w_{\mathcal{H}}({v}^3_k) = s_{1}+\sum_{j}a_{j} ,\\
		w_{\mathcal{H}}({v}^0_k) &=& w_{\mathcal{H}}({v}^2_k) = s_{2}+\sum_{j}b_{j},\\
	\end{eqnarray*}
	where $j$ varies from $1$ to $n-1,$ $j\neq k$ and $j\neq (k+\frac{n-1}{2})\mod (n-1)$.\\
	For all $j\neq 0$, by comparing $w_{\mathcal{H}}(v_{j}^{i})$ and $w_{\mathcal{H}}(v_{j+1}^{i})$, one can obtain,
	\begin{eqnarray*}
		a_{(j+1)~\text{mod}~(n-1)}+a_{(j+\frac{(n-1)}{2}+1)~\text{mod}~(n-1)}&=&a_{j}+a_{(j+\frac{n-1}{2})~\text{mod}~(n-1)} \\
		b_{(j+1)~\text{mod}~(n-1)}+b_{(j+\frac{(n-1)}{2}+1)~\text{mod}~(n-1)}&=&b_{j}+b_{(j+\frac{n-1}{2})~\text{mod}~(n-1)}.
	\end{eqnarray*}
	Now let,
	\begin{eqnarray*}
		a_{k}+a_{k+\bigl(\frac{n-1}{2}\bigr)} &=& a ~\textnormal{(say)},\\
		b_{k}+b_{k+\bigl(\frac{n-1}{2}\bigr)} &=& b~\textnormal{(say)},
	\end{eqnarray*}
	where $k \in \{1,2,...,\frac{n-1}{2}\}$.
	\par Since $c = w_{\mathcal{H}}(v_0^0) = \bigl( \frac{n-1}{2}\bigr) b = w_{\mathcal{H}}(v_0^1)=\bigl(\frac{n-1}{2}\bigr)a$, we have $a = b = a_{0}$. This leads, $s_{1}=s_{2}=s.$
	But the total weight is, $2s+(n-1)a_{0} = 2n(4n+1)$. Further, $c = \bigl(\frac{n-1}{2}\bigr)a_{0} = s+ \bigl(\frac{n-3}{2}\bigr)a_{0}$ and therefore $a_{0} = \frac{2n(4n+1)}{n+1} = 8n-6+\frac{6}{n+1}$, which is not an integer except for $n=5$.
	\vskip 0.5 cm
	\par Conversely, for $n=5$, Define,
	\begin{table}[htbp]
		\centering
		\begin{tabular}{cccc}
			$f(v_0^0)$ = 16, & $f(v_0^1)$ = 15, & $f(v_0^2)$ = 19, & $f(v_0^3)$ = 20, \\
			$f(v_1^0)$ = 18, & $f(v_1^1)$ = 6, & $f(v_2^2)$ = 1, & $f(v_3^3)$ = 7, \\
			$f(v_2^0)$ = 11, & $f(v_2^1)$ = 2, & $f(v_2^2)$ = 9, & $f(v_2^3)$ = 17, \\
			$f(v_3^0)$ = 3, & $f(v_3^1)$ = 14, & $f(v_3^2)$ = 13, & $f(v_3^3)$ = 8, \\
			$f(v_4^0)$ = 10, & $f(v_4^1)$ = 12, & $f(v_4^2)$ = 5, & $f(v_4^3)$ = 4. \\
		\end{tabular}
	\end{table}
	\begin{figure}[htbp]
		\centering
		\includegraphics[width=100mm]{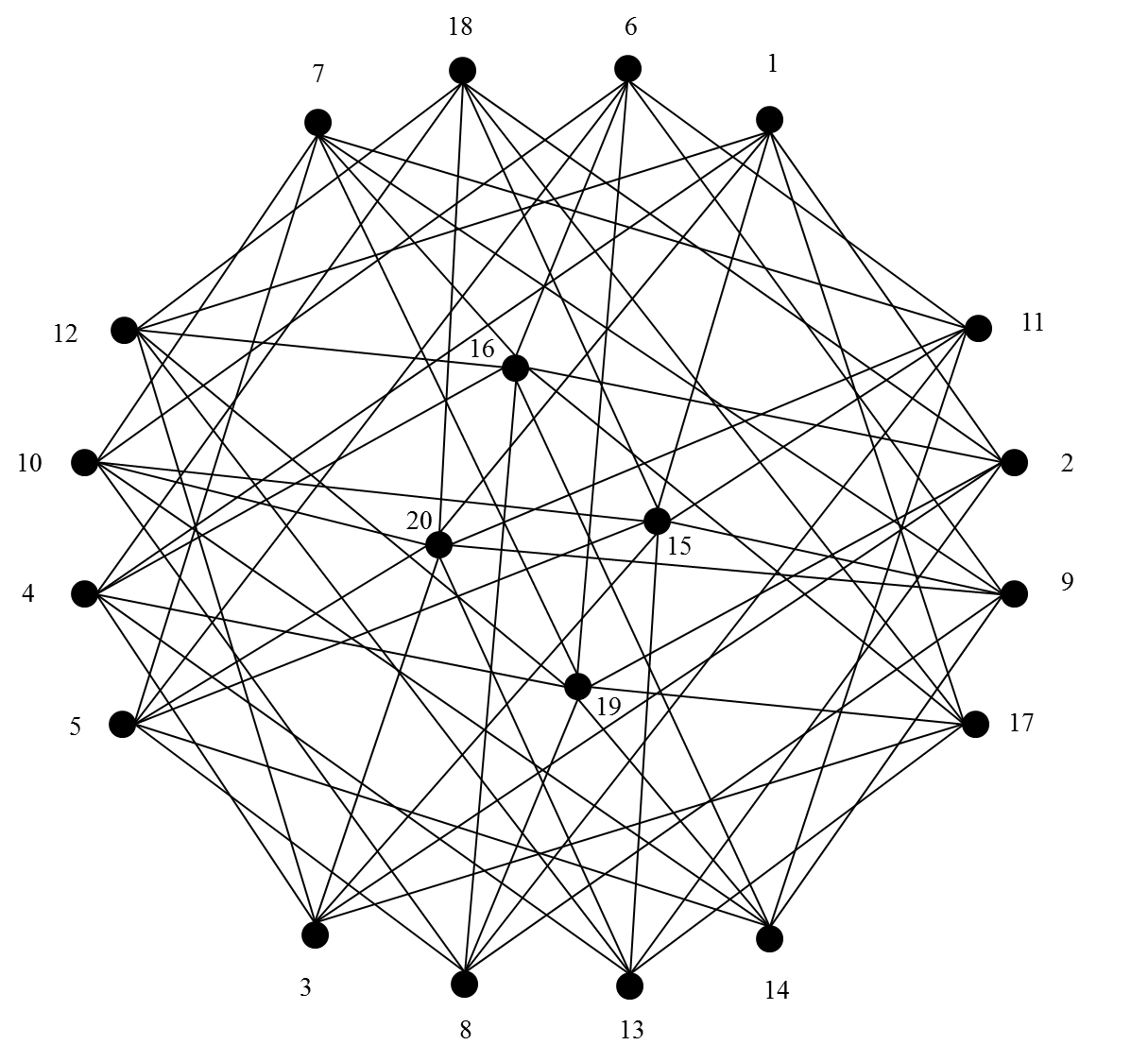}
		\caption{Distance magic labeling of $H_{3,5}\times C_{4}$}
		\label{fig1}
	\end{figure}
	\par Here $f$ is a distance magic labeling of $H_{3,5}\times C_{4}$, with magic constant 70.
\end{proof}
\begin{corollary}
	$H_{2n+1,2n+3} \times C_{4}$ is distance magic if and only if $n =5$.
\end{corollary}
\par From Theorem \ref{prbsol1} and \ref{prbsol2}, we observe that the distance magicness of the given graph is not sufficient for the above two problems to hold.
\section{$(a,d)$-distance antimagic labeling of $H_{m,n}$ }
\noindent This section exhibits the $(a,d)$-distance antimagic labeling of some subclasses of $H_{m,n}$ and provides a method for the construction of a non-regular $(a,d)$-distance antimagic graph from an existing one.
\begin{lemma}\label{lem2}
	Let $G$ be an $r$-regular graph on $n$ vertices.
	\begin{itemize}
		\item[\textnormal{(i)}] If $r$ is odd and $d$ is even, then $G$ is not $(a,d)$-distance antimagic.
		\item[\textnormal{(ii)}] If $r, n$ are even and $d$ is odd, then $G$ is not $(a,d)$-distance antimagic.
	\end{itemize}
\end{lemma}
\begin{proof}
	If $G$ is an $r$-regular $(a,d)$-distance antimagic graph, then by using Lemma \ref{siman1},
	\begin{equation}\label{eqn1}
	a =  \frac{n(r-d)+(r+d)}{2}
	\end{equation}
	Right-hand side of \ref{eqn1} becomes an integer only when $r$ and $d$ are of same parity.
\end{proof}
\par From Lemma \ref{lem2}, it is observed that when $m$ and $n$ are even, $H_{m,n}$ is not $(a,d)$-distance antimagic for any odd $d$. On the other hand, if $m$ is odd and $n$ is even, $H_{m,n}$ is not $(a,d)$-distance antimagic for any even $d$. \par We know that $K_{n}$ is an $(a,d)$-distance antimagic graph if and only if $d=1$. Moreover, if $M$ is any perfect matching of $K_{n}$, then any $(n-2)$-regular graph $G$ on $n$ vertices is isomorphic to $K_{n}-M$. Since the twin vertices of $G$ share a common neighborhood, $G$ is not $(a,d)$-distance antimagic for any $d\geq1.$ Now, for even $n$, the following lemma gives a necessary condition for an $(n-3)$-regular graph on $n$ vertices to be $(a,d)$-distance antimagic.
\begin{lemma}\label{lem3}
	Let $n\geq4$ be an even integer. If $G$ is an $(n-3)$-regular $(a,d)$-distance antimagic graph on $n$ vertices, then $d=1.$
\end{lemma}
\begin{proof}
	From Lemma \ref{siman1}, we have,
	\begin{equation*}
	d\leq \frac{3(n-3)}{n-1} < 3.
	\end{equation*}
	Again, by Lemma \ref{lem2}, $d=1$ is the only possibility.
\end{proof}
\begin{theorem}
	$H_{2n,2n+3}$ is $(2n^2+3n-1,1)$-distance antimagic.
\end{theorem}
\begin{proof}
	Let $G\cong H_{2n,2n+3}$. Consider a function $f$ on $V(G)$ as,
	\[  f(v_{i}) = \left\{
	\begin{array}{ll}
	2n+3 - i  & $ for $ 0\leq i \leq n+1 $ and $ i = 2n+2\\
	i-n  & $ for $ n+2 \leq i \leq 2n+1,  \\
	\end{array}
	\right.\]
	for which the distinct weights in arithmetic progression are,
	{ \[  w_{G}(v_{i}) = \left\{
		\begin{array}{ll}
		2n^{2}+2n-3+i  & $ for $ n+2 \leq i \leq 2n+1 \\
		2n^2+4n-1+\beta_{0}(i) & $ for $ i = 0,n+1 \\
		2n^{2}+5n-i & $ for $ 1 \leq i \leq n-1  \\
		2n^{2}+5n+1-\beta_{n}(i)  & $ for $ i = 2n+2,n. \mbox{~~ }
		\end{array}
		\right.\]}	
\end{proof}
\begin{theorem}
	$H_{4n+1,4n+4}$ is $(a,d)$-distance antimagic if and only if $d = 1$.
\end{theorem}
\begin{proof}
	Let $G\cong H_{4n+1,4n+4}.$ Consider a function $f$ on $V(G)$ as,
	\[  f(v_{i}) = \left\{
	\begin{array}{ll}
	4n+4 - i  & $ for $ 0\leq i \leq 2n+1 \\
	i-2n-2\alpha(i)  & $ for $ 2n+2 \leq i \leq 4n+3. \\
	\end{array}
	\right.\]
	Notice that the distinct weights in arithmetic progression are,
	{ \small	\[  w_{G}(v_{i}) = \left\{
		\begin{array}{ll}
		8n^{2}+10n+1  & $ for $ i = 2n+1 \\
		8n^{2}+10n+2  & $ for $ i = 2n+4 \\
		8n^{2}+8n-2\alpha(i+1)+i  & $ for $ {2n+3 \leq i \leq 4n+2},i\neq 2n+4,4n+1 \\
		8n^{2}+12n+1  & $ for $ i = 4n+1 \\
		8n^{2}+12n+2  & $ for $ i = 0 \\
		8n^{2}+12n+3  & $ for $ i = 2n-1 \\
		8n^{2}+12n+4  & $ for $ i = 2n+2\\
		8n^{2}+14n+2+4\alpha(i+1)-i & $ for $ 1 \leq i \leq 2n , i\neq 2,4,2n-1 \\
		8n^{2}+14n+2  & $ for $ i = 4\\
		8n^{2}+14n+3  & $ for $ i = 4n+3\\
		8n^{2}+14n+4 & $ for $ i = 2.
		\end{array}
		\right.\] }\\
	\par Hence, $f$ is an $(a,1)$-distance antimagic labeling with $a=8n^2+10n+1$.
	The converse follows from Lemma \ref{lem3}.
\end{proof}
\par It is observed from \cite{arumugam2016distance} that when $n$ is even, any $(n-2)$-regular distance magic graph on $n$ vertices can be extended to a non-regular distance magic graph on $n+1$ vertices. We know that if $G$ is $(a,d)$-distance antimagic graph then $G$ is trivially a distance antimagic graph. The following theorem provides a technique to construct new class of non-regular $(a,d)$-distance antimagic graphs from existing ones and also a partial solution to a much stronger version of the Problem \ref{prob3}.
\begin{theorem}\label{constthm2}
	If $G$ is an $(n-k)$-regular $(a,k-2)$-distance antimagic graph on $n$ vertices with $2 \leq k \leq \frac{1+\sqrt {8n-7}}{2},$ then there exists a non-regular $(a',k-2)$-distance antimagic graph on $n+1$ vertices, where $a' = a+n+1.$
\end{theorem}
\begin{proof}
	Let $G$ be an $(n-k)$-regular graph on $n$ vertices and $f$ be an $(a,k-2)$-distance antimagic labeling of $G$. Clearly, $a = \frac{n^2-2kn+3n-2}{2}$. If $k = 2$, the result follows from Theorem \ref{aru1}. When $k> 2$, construct a graph $G^{\dagger} \cong G+K_{1}$, where $u$ is the new vertex that induces $K_{1}$ of $G^{\dagger}$. Consider a function $f^{\dagger}$ as,
	\[  f^{\dagger}(v) = \left\{
	\begin{array}{ll}
	f(v) & $ for $ v \in V(G) \\
	n+1  & $ for $ v = u.  \\
	\end{array}
	\right.\]
	For every $v$ of $G^{\dagger}$ except $u$, $w_{G^{\dagger}}(v)= w_{G}(v)+ n+1$, which forms an arithmetic progression with $n^{th}$ term, $\frac{n^2+n}{2}-(k-2)$. Further, $w_{G^{\dagger}}(u)$ = $1+2+...+n$ = $\frac{n^2+n}{2}$. Therefore, $f^{\dagger}$ is a $(a',k-2)$-distance antimagic labeling of $G^{\dagger}$ with $a'=a+n+1.$
\end{proof}
\begin{corollary}\label{cor1}
	$H_{2n,2n+3} +K_{1}$ is $(2n^2+5n+3,1)$-distance antimagic.
\end{corollary}
\begin{corollary}\label{cor2}
	$H_{4n,4n+4} +K_{1}$ is $(8n^2+14n+6,1)$-distance antimagic.
\end{corollary}
\par It is interesting to see that Theorem \ref{constthm2}, Corollary \ref{cor1} and \ref{cor2} together exhibit $(a,1)$-distance antimagic graphs $G$ on $n$ vertices, for any even integer $n\geq6$ and for any odd integer $n\geq9$ and $n\equiv 1 ~\text{mod}~4$.
\par From Lemma \ref{lem3}, $H_{4,n}$ is not $(a,1)$-distance antimagic, when $n$ is even. The following theorem guarantees the existence of $(a,1)$-distance antimagic labeling of $H_{4,n}$, when $n$ is not even.
\begin{theorem}
	$H_{4,n}$ is $(a,1)$-distance antimagic if $n\equiv3 ~\textnormal{mod}~4.$
\end{theorem}
\begin{proof}
	Let $G\cong H_{4,n}$ with $n\equiv3 ~\textnormal{mod}~4$. 
	Define a temporary labeling $f$ on $G$ as,
	\[  f(v_{i}) = \left\{
	\begin{array}{ll}
	n  & $ for $ i=n-4 \\
	ki+1~\text{mod}~n  & $ for $ i\neq n-4,\\
	\end{array}
	\right.\]
	where, $k = \frac{1}{4}(n+1)$.
	Since $\gcd(k,n)=1$, $f$ is a bijection. The weights of the vertices of $G$ can be obtained as,
	\[  w_{G}(v_{i}) = \left\{
	\begin{array}{ll}
	2n+4+i  & $ for $ i\equiv 0 ~\text{mod}~4, 0 \leq i \leq n-7  \\
	2n+4+i  & $ for $ i\equiv 1 ~\text{mod}~4, 0 \leq i \leq n-10 \\
	n+4+i   & $ for $ i\equiv 2 ~\text{mod}~4, 0 \leq i \leq n-9 \\
	n+4+i   & $ for $ i\equiv 3 ~\text{mod}~4, 0 \leq i \leq n-8 \\
	n+3 & $ for $ i= n-1\\
	3n+1 & $ for $ i=n-3 \\
	2n-1 & $ for $ i=n-5\\
	2n-n\beta_{n-5}(i)+4+i & $ for $ i = n-6,n-4,n-2.
	\end{array}
	\right.\]

	\normalsize \par Define a new function $f^{\dagger}$ on $G$ as,
	\[  f^{\dagger}(v_{i}) = \left\{
	\begin{array}{ll}
	n- f(v_{i}) & $ for $ i\equiv0,3~\text{mod}~4, 0\leq i \leq n-7 $ and $ i=n-1 \\
	f(v_{i})  & $ otherwise$.\\
	\end{array}
	\right.\]
	Hence we have,
	\footnotesize	{\[  w^{\dagger}_{G}(v_{i}) = \left\{
		\begin{array}{ll}
		n+w_{G}(v_{i})-2f(v_{i-1}) & $ for $ i\equiv 0 ~\text{mod}~4, 0\leq i\leq n-7 \\
		3n+w_{G}(v_{i})-2f(v_{i-1})-2f(v_{i-2})-2f(v_{i+2})& $ for $ i\equiv 1 ~\text{mod}~4, 1\leq i\leq n-10 \\
		3n+w_{G}(v_{i})-2f(v_{i+1})-2f(v_{i+2})-2f(v_{i-2})& $ for $ i\equiv 2 ~\text{mod}~4, 2\leq i\leq n-9 \\
		n+w_{G}(v_{i})-2f(v_{i+1}) & $ for $ i\equiv 3 ~\text{mod}~4, 3\leq i\leq n-8 \\
		2n+1 & $ for $ i = n-1\\
		\frac{7}{2}(n-1)-(1-\beta_{n-5}(i))-i & $ for $ i = n-2,n-3,n-5,n-6\\
		2n & $ for $ i = n-4.\\
		\end{array}
		\right.\]}
	\normalsize
	\par After simplification,
	\[  w^{\dagger}_{G}(v_{i}) = \left\{
	\begin{array}{ll}
	\frac{3n+5+i}{2} & $ for $ i\equiv 0 ~\text{mod}~4, 0\leq i\leq n-7 \\
	\frac{4n-3-i}{2} & $ for $ i\equiv 1 ~\text{mod}~4, 1\leq i\leq n-10 \\
	\frac{5n-5-i}{2} & $ for $ i\equiv 2 ~\text{mod}~4, 2 \leq i\leq n-9 \\
	\frac{4n+3+i}{2} & $ for $ i\equiv 3 ~\text{mod}~4, 3 \leq i\leq n-8 \\
	2n+1 & $ for $ i = n-1\\
	\frac{7}{2}(n-1)-(1-\beta_{n-5}(i))-i & $ for $ i = n-2,n-3,n-5,n-6\\
	2n & $ for $ i = n-4.\\
	\end{array}
	\right.\]
	Here, $f^{\dagger}$ is an $(a,1)$-distance antimagic labeling of $G$ with $a = \frac{3n+5}{2}.$
\end{proof}
\section{Conclusion and Scope}
\noindent This paper discusses the concept of distance magicness and $(a,d)$-distance antimagicness of subfamilies of Harary graph. Two techniques for building larger classes of distance magic and $(a,d)$-distance antimagic graphs from the existing ones, are constructed. More importantly, some classes of non-regular distance magic and $(a,d)$-distance antimagic graphs are constructed through these techniques.
\par A comprehensive list of subfamilies of $H_{m,n}$, which are known to admit distance magic and $(a,d)$-distance antimagic labeling, is given in Table 1.
\begin{table}[htbp]
	\scriptsize
	\centering
	\caption{ Distance magic and $(a,d)$-distance antimagic labeling of $H_{m,n}.$}
	\label{tab:1}
	\begin{tabular}{lll}
		\hline\noalign{\smallskip}
		& Distance magic & $(a,d)$-distance \\ & & antimagic  \\
		\noalign{\smallskip}\hline\noalign{\smallskip}
		$H_{2,n}$ & Yes~(only if $n=4$)& Yes~(only if $n$-odd, $d=1$)\\
		$H_{3,n}$ & Yes~(only if $n=5$) & \\
		$H_{4,n}$ & Yes~(only if $n=6)$  &  Yes~(if $n\equiv3~\text{mod}~4$, $d=1$), No~(if $n$ is even)  \\
		$H_{6,n}$ & Yes~(only if $n=8,24)$  &  No~(if $n$-even, $d$-odd)  \\
		$H_{10,n}$ & Yes~(only if $n=12,20,60)$  &  No~(if $n$-even, $d$-odd)  \\
		$H_{2n,2n+2} $& Yes  &  No~(if $d$-odd)  \\
		$H_{2n+1,2n+3}$ & Yes  & No \\
		$H_{2n,2n+3}$ &   & Yes~($d=1$) \\
		$ H_{n,n+1}$ & No  &  Yes~(only if $d=1$)  \\
		$H_{4n+1,4n+4}$ &  & Yes~(only if $d=1$) \\
		$ H_{m,n},$($m$-odd, $n$-even) & No & No~(if $d$-even) \\
		$H_{2m,n}$, where & & \\
		$\text{(i)}$ $m$-odd,\\ $2m^{2}+m \equiv 0 ~\text{mod}~n$,\\ $\frac{n}{gcd(n,m+1)} \equiv 0 ~\text{mod}~2$ \\ and $n \geq 2m+2 $& Yes  & No(if $d$-even) \\
		$\textnormal(ii)$ $\frac{n}{gcd(n,m+1)} \equiv 1 ~\text{mod}~2$ & No &  \\
		$\textnormal(iii)$ $m$-odd, $2m^{2}+m\not\equiv 0 ~\text{mod}~n$ & No &  \\
		\noalign{\smallskip}\hline
	\end{tabular}
\end{table}


\begin{thebibliography}{99}





\bibitem{bondy1976graph}
{\sc J. A. Bondy} and {\sc U. S. R. Murty}, {Graph Theory with Applications}, American Elsevier Publishing Co., Inc., New York, 1976.


\bibitem{hammack2011handbook}
{\sc R. Hammack, W. Imrich} and {\sc S. Klav{\v{z}}ar}, Handbook of Product Graphs, CRC Press, Boca Raton, FL, 2011.



\bibitem{vilfred}
{\sc V. Vilfred}, $\sum-$ labelled graphs and circulant graphs, Ph.D. thesis, University
of Kerala, Trivandrum, India, 1994.

\bibitem{miller2003distance}
{\sc M. Miller, C. Rodger} and {\sc R. Simanjuntak}, Distance magic labelings of
graphs, {\it Australasian Journal of Combinatorics} {\bf 28} (2003), 305--315.

\bibitem{sugeng2009distance}
{\sc K. Sugeng, D. Froncek, M. Miller, T. Ryan} and {\sc J. Walker},  On distance
magic labeling of graphs, {\it Journal of Combinatorial Mathematics and Combinatorial Computing} {\bf 71} (2009), 39--48.

\bibitem{froncek2006fair}
{\sc D. Froncek, P. Kov{\'a}r} and {\sc T. Kov{\'a}rov{\'a}}, Fair incomplete tournaments, {\it Bulletin of the Institute of Combinatorics and its Applications} {\bf 48} (2006), 31--33.

\bibitem{froncek2007fair}
{\sc D. Froncek}, Fair incomplete tournaments with odd number of teams and large number of games, {\it Congressus Numerantium} {\bf 187} (2007), 83--89.

\bibitem{jinnah}
{\sc M. I. Jinnah}, On $\sum-$ labelled graphs, in: Technical Proceedings of
Group Discussion on Graph Labeling Problems, (eds.) B.D. Acharya
and S.M. Hedge (1999), 71--77.

\bibitem{rao}
{\sc S. B. Rao}, Sigma graphs-a survey, in: Labelings of Discrete Structures
and Applications, (eds.) B.D. Acharya, S. Arumugam, A. Rosa, Narosa
Publishing House, New Delhi (2008) 135--140.

\bibitem{anholcer2015distance}
{\sc M. Anholcer, S. Cichacz, I. Peterin} and {\sc A. Tepeh}, Distance magic labeling and two products of graphs, {\it Graphs and Combinatorics} {\bf 31} (2015), 1125--1136.

\bibitem{cichacz2013distance}
{\sc S. Cichacz},  Distance magic $(r,t)$ hypercycles, {\it Utilitas Mathematics} {\bf 101} (2016), 283--294.

\bibitem{shafiq2009distance}
{\sc M. K. Shafiq, G. Ali} and {\sc R. Simanjuntak}, Distance magic labelings of a union of graphs, {\it AKCE International Journal of Graphs and Combinatorics} {\bf 6} (2009), 191--200.

\bibitem{cichacz2018constant}
{\sc S. Cichacz} and {\sc A. G{\H{o}}rlich},  Constant sum partition of sets of integers and distance magic graphs, {\it Discussiones Mathematicae Graph Theory} {\bf 38} (2018), 97--106.

\bibitem{arumugam2016distance}
{\sc S. Arumugam, N. Kamatchi} and {\sc P. Kov{\'a}r}, Distance magic graphs, {\it Utilitas Mathematica} {\bf 99} (2016), 131--142.


\bibitem{arumugam2012d}
{\sc S. Arumugam} and {\sc N. Kamatchi}, On (a,d)-distance antimagic graphs, {\it Australasian Journal of Combinatorics} {\bf 54} (2012), 279--287.

\bibitem{kamatchi2013distance}
{\sc N. Kamatchi} and {\sc S. Arumugam}, Distance antimagic graphs, {\it Journal of Combinatorial Mathematics and Combinatorial Computing} {\bf 84} (2013) 61--67.

\bibitem{simanjuntak2013distance}
{\sc R. Simanjuntak} and {\sc K. Wijaya}, On distance antimagic graphs, {\it arXiv preprint arXiv: 1312.7405} (2013).


\bibitem{kovar2012note}
{\sc P. Kov{\'a}r, D. Froncek} and {\sc T. Kov{\'a}rov{\'a}}, A note on 4-regular distance magic graphs, {\it Australasian Journal of Combinatorics} {\bf 54} (2012) 127--132.












							
							
							

							

							

								
								
								
										
%
%
%
%
%
%
%
%
%
%
%
%
%
%
%
%
%
%
%
%
%
%
%
%
%
%
%
%
%
%
%
%
%
%
%
%
%
%
%
%
%
%


\end{thebibliography}
\end{document}